\documentclass[runningheads]{llncs}
\usepackage[linesnumberedhidden,lined,ruled]{algorithm2e}

\usepackage[colorlinks=true,urlcolor=blue,linkcolor=blue,citecolor=blue]{hyperref}
\usepackage{amsmath}
\usepackage{amssymb}
\usepackage{doi}
\usepackage{tikz,tikz-3dplot}
\usetikzlibrary{graphs}
\usetikzlibrary{graphs.standard}
\newcommand{\R}{\mathbb{R}}

\DeclareMathOperator\conv {conv} 
\DeclareMathOperator\TSP {TSP} 
\DeclareMathOperator\ATSP {ATSP}
\usepackage{graphicx}
\begin{document}
\title{Simulated annealing approach to verify vertex adjacencies in the traveling salesperson polytope}
\titlerunning{SA approach to verify vertex adjacencies in the TSP polytope}
%
\author{Anna Kozlova\inst{1} \and
Andrei Nikolaev\inst{1}}
\authorrunning{Kozlova A.\,P, Nikolaev A.\,V.}
%
\institute{P.\,G. Demidov Yaroslavl State University\\
\email{fyz95@mail.ru, andrei.v.nikolaev@gmail.com}}
\maketitle              
\begin{abstract}
We consider 1-skeletons of the symmetric and asymmetric traveling salesperson polytopes whose vertices are all possible Hamiltonian tours in the complete directed or undirected graph, and the edges are geometric edges or one-dimensional faces of the polytope.
It is known that the question whether two vertices of the symmetric or asymmetric traveling salesperson polytopes are nonadjacent is NP-complete.
A sufficient condition for nonadjacency can be formulated as a combinatorial problem: if from the edges of two Hamiltonian tours we can construct two complementary Hamiltonian tours, then the corresponding vertices of the traveling salesperson polytope are not adjacent. 
We consider a heuristic simulated annealing approach to solve this problem. It is based on finding a vertex-disjoint cycle cover and a perfect matching.
The algorithm has a one-sided error: the answer ``not adjacent'' is always correct, and was tested on random and pyramidal Hamiltonian tours.

\end{abstract}

\section{Introduction} 

We consider a classical traveling salesperson problem on a complete directed or undirected graph. 

\textsc{Symmetric traveling salesperson problem.} Given a complete weighted graph $K_n=(V,E)$, it is required to find a Hamiltonian cycle of minimum weight.

\textsc{Asymmetric traveling salesperson problem.} Given a complete weighted digraph $D_n=(V,A)$, it is required to find a Hamiltonian tour of minimum weight.


We denote by $HC_n$ the set of all Hamiltonian cycles in $K_n$ and by $HT_n$ the set of all Hamiltonian tours in $D_n$. With each Hamiltonian cycle $x \in HC_n$ we associate a characteristic vector $x^v \in \R^{E}$ by the following rule:
\[
x^v_e = 
\begin{cases}
1,& \text{ if the cycle } x \text{ contains an edge } e \in E,\\
0,& \text{ otherwise. }
\end{cases}
\]
With each Hamiltonian tour $y \in HT_n$ we associate a characteristic vector $y^v \in \R^{A}$ by the following rule:
\[
y^v_a = 
\begin{cases}
1,& \text{ if the tour } y \text{ contains an edge } a \in A,\\
0,& \text{ otherwise. }
\end{cases}
\]
The polytope 
\[\TSP(n) = \conv \{x^v \ | \ x \in HC_n \}\]
is called \textit{the symmetric traveling salesperson polytope},  
and the polytope 
\[\ATSP(n) = \conv \{y^v \ | \ y \in HT_n \}\]
is called \textit{the asymmetric traveling salesperson polytope}.

The $1$-skeleton of a polytope $P$ is the graph whose vertex set is the vertex set of $P$ (characteristic vectors $x^v$ for the traveling salesperson problem) and edge set is the set of geometric edges or one-dimensional faces of $P$. Many papers are devoted to the study of $1$-skeletons associated with combinatorial problems. On the one hand, the vertex adjacency in $1$-skeleton is of great interest for the development of algorithms to solve problems based on local search technique (when we choose the next solution as the best one among adjacent solutions). For example, various algorithms for perfect matching, set covering, independent set, ranking of objects, problems with fuzzy measures, and many others are based on this idea \cite{Aguilera,Balinski,Chegireddy,Combarro,Gabow,Matsui}. On the other hand, some characteristics of $1$-skeletons, such as the diameter and the clique number, estimate the time complexity for different computation models and classes of algorithms \cite{Bondarenko,Bondarenko-Maksimenko,Bondarenko-Nikolaev-2017,Grotschel}.

Unfortunately, the classical result by Papadimitriou states that the construction of 1-skeleton of the traveling salesperson polytope is NP-complete for both directed and undirected graphs.
\begin{theorem} [Papadimitriou, \cite{Papadimitriou}]
    The question whether two vertices of the polytopes $\TSP(n)$ or $\ATSP(n)$ are nonadjacent is NP-complete. 
\end{theorem}

However, the vertex adjacency test for $\TSP(n)$ and $\ATSP(n)$ is an interesting problem itself. Note that a geometric approach to the construction of 1-skeleton of the traveling salesperson polytope seems not very promising because both polytopes have superexponential number of vertices and faces \cite{Fiorini,Grotschel}.

In \cite{Rao} the sufficient condition for vertex adjacency in the traveling salesperson polytope was reformulated in a combinatorial form.

\begin{lemma}[Sufficient condition for nonadjacency] \label{lemma_nonadjacency}
	If from the edges of \break
	two Hamiltonian tours $x$ and $y$ it is possible to construct two complementary Hamiltonian tours $z$ and $w$, then the corresponding vertices $x^v$ and $y^v$ of the polytope $\TSP(n)$ (or $\ATSP(n)$) are not adjacent. 
\end{lemma}

\begin{figure}[t]
	\centering
	\begin{tikzpicture}[scale=0.85]
	\begin{scope}[every node/.style={circle,thick,draw,inner sep=3pt}]
	\node (A) at (0,0) {1};
	\node (B) at (1,0) {2};
	\node (C) at (2,0) {3};
	\node (D) at (3,0) {4};
	\node (E) at (4,0) {5};
	\node (F) at (5,0) {6};
	\node (G) at (6,0) {7};
	\node (H) at (7,0) {8};
	\end{scope}
	\draw [thick] (A) edge (B);
	\draw [thick] (B) edge [bend left=50] (D);
	\draw [thick] (D) edge [bend left=40] (G);
	\draw [thick] (G) edge (F);
	\draw [thick] (F) edge [bend left=50] (H);
	\draw [thick] (H) edge [bend left=45] (E);
	\draw [thick] (E) edge [bend left=50] (C);
	\draw [thick] (C) edge [bend left=50] (A);
	\draw (-1, 0) node{\textit{x}};
	\end{tikzpicture}
	\\
	\begin{tikzpicture}[scale=0.85]
	\begin{scope}[every node/.style={circle,thick,draw,inner sep=3pt}]
	\node (A) at (0,0) {1};
	\node (B) at (1,0) {2};
	\node (C) at (2,0) {3};
	\node (D) at (3,0) {4};
	\node (E) at (4,0) {5};
	\node (F) at (5,0) {6};
	\node (G) at (6,0) {7};
	\node (H) at (7,0) {8};
	\end{scope}
	\draw [thick,dashed] (A) edge (B);
	\draw [thick,dashed] (B) edge (C);
	\draw [thick,dashed] (C) edge (D);
	\draw [thick,dashed] (D) edge [bend left=50] (F);
	\draw [thick,dashed] (F) edge (G);
	\draw [thick,dashed] (G) edge (H);
	\draw [thick,dashed] (H) edge [bend left=45] (E);
	\draw [thick,dashed] (E) edge [bend left=40] (A);
	\draw (-1, 0) node{\textit{y}};
	\end{tikzpicture}
	\\
	\begin{tikzpicture}[scale=0.85]
	\begin{scope}[every node/.style={circle,thick,draw,inner sep=3pt}]
	\node (A) at (0,0) {1};
	\node (B) at (1,0) {2};
	\node (C) at (2,0) {3};
	\node (D) at (3,0) {4};
	\node (E) at (4,0) {5};
	\node (F) at (5,0) {6};
	\node (G) at (6,0) {7};
	\node (H) at (7,0) {8};
	\end{scope}
	\draw [thick] (A) edge (B);
	\draw [thick] (B) edge [bend left=50] (D);
	\draw [thick,dashed] (D) edge [bend left=50] (F);
	\draw [thick,dashed] (F) edge (G);
	\draw [thick,dashed] (G) edge (H);
	\draw [thick,dashed] (H) edge [bend left=45] (E);
	\draw [thick] (E) edge [bend left=50] (C);
	\draw [thick] (C) edge [bend left=50] (A);
	\draw (-1, 0) node{\textit{z}};
	\end{tikzpicture}
	\\
	\begin{tikzpicture}[scale=0.85]
	\begin{scope}[every node/.style={circle,thick,draw,inner sep=3pt}]
	\node (A) at (0,0) {1};
	\node (B) at (1,0) {2};
	\node (C) at (2,0) {3};
	\node (D) at (3,0) {4};
	\node (E) at (4,0) {5};
	\node (F) at (5,0) {6};
	\node (G) at (6,0) {7};
	\node (H) at (7,0) {8};
	\end{scope}
	\draw [thick,dashed] (A) edge (B);
	\draw [thick,dashed] (B) edge (C);
	\draw [thick,dashed] (C) edge (D);
	\draw [thick] (D) edge [bend left=40] (G);
	\draw [thick] (G) edge (F);
	\draw [thick] (F) edge [bend left=50] (H);
	\draw [thick] (H) edge [bend left=45] (E);
	\draw [thick,dashed] (E) edge [bend left=40] (A);
	\draw (-1, 0) node{\textit{w}};
	\end{tikzpicture}
	\caption{Two complementary tours $z$ and $w$ are constructed from the edges of $x$ and $y$}
	\label{image1}
\end{figure}
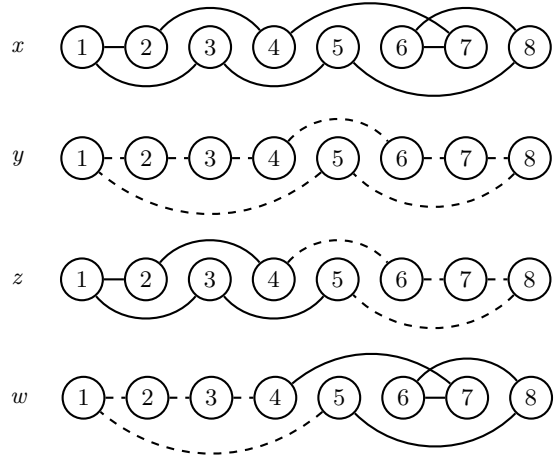

From the geometric point of view, the Lemma~\ref{lemma_nonadjacency} means that the segment connecting two vertices $x^v$ and $y^v$ intersects with the segment connecting two other vertices $z^v$ and $w^v$ of the polytope $\TSP(n)$ (or $\ATSP(n)$ correspondingly), thus it cannot be an edge in 1-skeleton. An example of a satisfied sufficient condition is shown in Fig.~\ref{image1}. 

Let us formulate the sufficient condition for vertex nonadjacency of the traveling salesperson polytope in the form of a combinatorial problem. \\ 
\textsc{Instance.} Let $x$ and $y$ be two Hamiltonian tours.\\ 
\textsc{Question.} Does the multigraph $x \cup y$ include a pair of Hamiltonian tours $z$ and $w$ different from $x$ and $y$ such that 
\[z \cup w = x \cup y \text{ and } z \cap w = \emptyset?\]

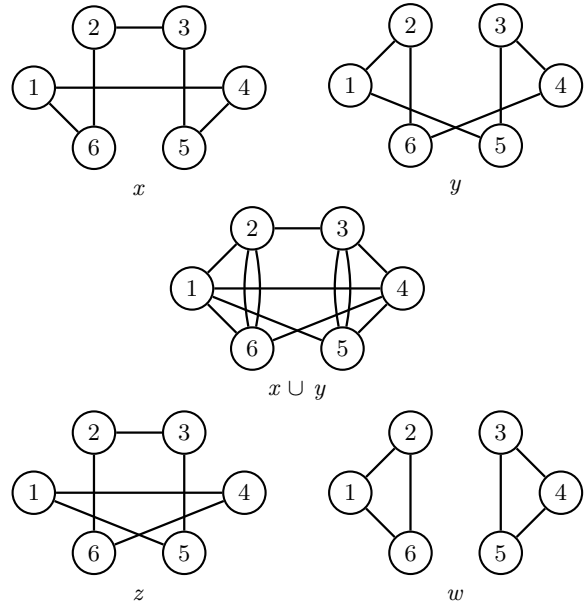
\begin{figure}[t]
	\centering
	\begin{tikzpicture}[scale=0.8]
	\begin{scope}[every node/.style={circle,thick,draw}]
	\node (A) at (0,0) {1};
	\node (B) at (1,1) {2};
	\node (C) at (2.5,1) {3};
	\node (D) at (3.5,0) {4};
	\node (E) at (2.5,-1) {5};
	\node (F) at (1,-1) {6};
	\end{scope}
	\draw [line width=0.3mm] (A) edge (D);
	\draw [line width=0.3mm] (D) edge (E);
	\draw [line width=0.3mm] (E) edge (C);
	\draw [line width=0.3mm] (C) edge (B);
	\draw [line width=0.3mm] (B) edge (F);
	\draw [line width=0.3mm] (F) edge (A);
	\draw (1.75, -1.7) node{\textit{x}};
	\end{tikzpicture}
	\hspace*{6mm}
	\begin{tikzpicture}[scale=0.8]
	\begin{scope}[every node/.style={circle,thick,draw}]
	\node (A) at (0,0) {1};
	\node (B) at (1,1) {2};
	\node (C) at (2.5,1) {3};
	\node (D) at (3.5,0) {4};
	\node (E) at (2.5,-1) {5};
	\node (F) at (1,-1) {6};
	\end{scope}
	\draw [line width=0.3mm] (A) edge (B);
	\draw [line width=0.3mm] (B) edge (F);
	\draw [line width=0.3mm] (F) edge (D);
	\draw [line width=0.3mm] (D) edge (C);
	\draw [line width=0.3mm] (C) edge (E);
	\draw [line width=0.3mm] (E) edge (A);
	\draw (1.75, -1.7) node{\textit{y}};
	\end{tikzpicture}
	\\
	\begin{tikzpicture}[scale=0.8]
	\begin{scope}[every node/.style={circle,thick,draw}]
	\node (A) at (0,0) {1};
	\node (B) at (1,1) {2};
	\node (C) at (2.5,1) {3};
	\node (D) at (3.5,0) {4};
	\node (E) at (2.5,-1) {5};
	\node (F) at (1,-1) {6};
	\end{scope}
	\draw [line width=0.3mm] (A) edge (D);
	\draw [line width=0.3mm] (D) edge (E);
	\draw [line width=0.3mm, bend right=10] (E) edge (C);
	\draw [line width=0.3mm] (C) edge (B);
	\draw [line width=0.3mm, bend right=10] (B) edge (F);
	\draw [line width=0.3mm] (F) edge (A);
	\draw [line width=0.3mm] (A) edge (B);
	\draw [line width=0.3mm, bend left=10] (B) edge (F);
	\draw [line width=0.3mm] (F) edge (D);
	\draw [line width=0.3mm] (D) edge (C);
	\draw [line width=0.3mm, bend right=10] (C) edge (E);
	\draw [line width=0.3mm] (E) edge (A);
	\draw (1.75, -1.7) node{\textit{x $\cup$ y}};
	\end{tikzpicture}
	\\
	\begin{tikzpicture}[scale=0.8]
	\begin{scope}[every node/.style={circle,thick,draw}]
	\node (A) at (0,0) {1};
	\node (B) at (1,1) {2};
	\node (C) at (2.5,1) {3};
	\node (D) at (3.5,0) {4};
	\node (E) at (2.5,-1) {5};
	\node (F) at (1,-1) {6};
	\end{scope}
	\draw [line width=0.3mm] (A) edge (D);
	\draw [line width=0.3mm] (D) edge (F);
	\draw [line width=0.3mm] (F) edge (B);
	\draw [line width=0.3mm] (B) edge (C);
	\draw [line width=0.3mm] (C) edge (E);
	\draw [line width=0.3mm] (E) edge (A);
	\draw (1.75, -1.7) node{\textit{z}};
	\end{tikzpicture}
	\hspace*{6mm}
	\begin{tikzpicture}[scale=0.8]
	\begin{scope}[every node/.style={circle,thick,draw}]
	\node (A) at (0,0) {1};
	\node (B) at (1,1) {2};
	\node (C) at (2.5,1) {3};
	\node (D) at (3.5,0) {4};
	\node (E) at (2.5,-1) {5};
	\node (F) at (1,-1) {6};
	\end{scope}
	\draw [line width=0.3mm] (A) edge (B);
	\draw [line width=0.3mm] (B) edge (F);
	\draw [line width=0.3mm] (F) edge (A);
	\draw [line width=0.3mm] (C) edge (D);
	\draw [line width=0.3mm] (D) edge (E);
	\draw [line width=0.3mm] (E) edge (C);
	\draw (1.75, -1.7) node{\textit{w}};
	\end{tikzpicture}
	\caption{An example of $w = (z \cup y) \backslash z$ that is not a Hamiltonian tour}
	\label{image2}
\end{figure}

By $x \cup y$ we denote a multigraph that contains all edges of both tours $x$ and $y$ (Fig.~\ref{image2}).

In this formulation, the problem is close to the $2$-peripatetic salesperson problem in which it is required to find two Hamiltonian tours of minimum weight without common edges. The $2$-peripatetic salesperson is NP-complete even for $4$-regular graphs \cite{DeKort}. Much attention was paid to the development of approximation algorithms for this problem (see, for example, \cite{Ageev,Baburin,Glebov}). 

However, the combinatorial form of the sufficient condition for nonadjacency has a number of differences from the $2$-peripatetic salesperson problem: 
\begin{itemize}
	\item this is a decision problem, not an optimization one; 
	\item the graph is a $4$-regular graph (or digraph) of a special form constructed as a union of two Hamiltonian tours; 
	\item it is required to find two Hamiltonian tours different from $x$ and $y$. 
\end{itemize}

Note also that if from the edges of two Hamiltonian tours $x$ and $y$ it is possible to construct another Hamiltonian tour $z$, then all the remaining edges $(x \cup y) \backslash z$ are almost certainly do not form a Hamiltonian tour (Fig.~\ref{image2}). Thus, instead of algorithms for a single Hamiltonian tour in the multigraph $x \cup y$, in this paper we consider a heuristic simulated annealing approach to test vertex adjacencies in the symmetric and asymmetric traveling salesperson polytopes based on finding a vertex-disjoint cycle cover and a perfect matching. 

\section{Simulated annealing} 
The simulated annealing borrows the concept from annealing in metallurgy where a metal material is repeatedly heated, kneaded, and cooled to enlarge the size of its crystals to
eliminate defects~\cite{Kirkpatrick}.

We consider a general scheme of the Algorithm~\ref{simmanneal}.
It is required to minimize the energy function specified for the current system state. The algorithm starts from a certain initial state: at each step a neighbor candidate state is generated which energy is compared to the energy of the previous state. If the energy decreases, the system transits to the new state, otherwise it may transit with a certain probability (to prevent falling into the local minimum). 

The algorithm receives input data in one of the following formats: 
\begin{enumerate}
    \item Two Hamiltonian tours $x = [a_1, \ldots, a_N]$ and $y = [b_1, \ldots, b_N]$, given as the permutations of vertices in a complete graph (or digraph) $K_N$; 
    \item $2/4$-regular graph ($2$~--- for directed and $4$~--- for undirected graphs) of size $N$, i.e. the union of two Hamiltonian tours, given as the adjacency list. 
\end{enumerate}

Other input parameters: the initial value of temperature $initT$, the maximum number of iterations $iterN$, and the the size of a queue of fixed edges $fixEdgesN$.
The algorithm stops when the solution is found or when the number of iterations exceeds the value of the parameter $iterN$.
As an output, the algorithms returns two complementary Hamiltonian tours $z$ and $w$ constructed from the edges of $x$ and $y$. By the sufficient condition (Lemma~\ref{lemma_nonadjacency}), the corresponding vertices $x^v$ and $y^v$ of the traveling salesperson polytope are not adjacent. 
If the algorithm cannot find the complementary tours, then it returns that the corresponding vertices are probably adjacent. Thus, the algorithm has a one-sided error: the answer ``not adjacent'' is always correct, while the answer ``probably adjacent'' leaves the possibility that the vertices actually are not adjacent. 

\SetAlgorithmName{Algorithm}{Algorithm list}{} 
\SetAlgoCaptionSeparator{.} 
\DontPrintSemicolon
\SetKwProg{Proc}{Procedure}{}{End}
\SetKw{Return}{Return}
\SetKwFor{For}{For}{}{End}
\SetKw{KwTo}{to}
\SetKw{KwOr}{or}
\SetKwIF{If}{ElseIf}{Else}{If}{Then}{Else}{Else}{End}
\begin{algorithm}[t]
    \label{simmanneal}
	\caption{Simulated Annealing Algorithm}
	\SetKwInOut{Input}{Input}
	\SetKwInOut{Output}{Output}
	\SetKwData{x}{$x$}
	\SetKwData{y}{$y$}
	\SetKwData{combinedGraph}{$combG$}
	\SetKwData{startTemp}{$initT$}
	\SetKwData{finishTemp}{$finalT$}
	\SetKwData{numOfIterations}{$iterN$}
	\SetKwData{numOfEdgesToFix}{$fixEdgesN$}
	\SetKwData{T}{$T$}
	\SetKwData{candidateEnergy}{$candE$}
	\SetKwData{zCand}{$zCand$}
	\SetKwData{wCand}{$wCand$}
	\SetKwData{currentEnergy}{$currE$}
	\SetKwData{z}{$z$}
	\SetKwData{w}{$w$}
	\SetKwData{k}{$k$}
	\SetKwData{true}{true}
	\SetKwData{false}{false}
	\SetKwData{or}{$OR$}
	\SetKwData{and}{$AND$}
	\SetKwData{SimulatedAnnealing}{$SimulatedAnnealing$}
	\SetKwData{getInitialState}{$GetInitialState$}
	\SetKwData{calculateEnergy}{$CalculateEnergy$}
	\SetKwData{transitToNewState}{$TransitToNewState$}
	\SetKwData{generateNeighbourCandidate}{$GenerateNeighbourCandidate$}
	\SetKwFunction{shouldAcceptCandidate}{$ShouldAcceptCandidate$}
	\SetKwFunction{coolingSchedule}{$CoolingSchedule$}
	\SetKwData{TestVertexAdjacency}{$TestVertexAdjacency$}
	
	\Input{ 
	    Hamiltonian tours \x and \y (or $2$/$4$-regular graph \combinedGraph), \newline 
	    initial temperature \startTemp, 
	    number of iterations \numOfIterations, size of a queue of fixed edges \numOfEdgesToFix
	}
	\Output{ vertices $x^v$ and $y^v$ are adjacent or not adjacent, complementary Hamiltonian tours \z and \w, if exist}
	\BlankLine
	\Proc{\SimulatedAnnealing{\x, \y, \combinedGraph, \startTemp, 
	\numOfEdgesToFix}}{
	    \T $\leftarrow$ \startTemp\;
		\z, \w $\leftarrow$ \getInitialState{\x, \y, \combinedGraph}\;
			
		\For{$k\leftarrow 1$ \KwTo \numOfIterations}{
    		    \If{\z and \w are Hamiltonian tours different from \x and \y}{
    		        \Return \z and \w\;
    		    }
    		  
		        \zCand, \wCand $\leftarrow$ \generateNeighbourCandidate{\z, \w, \numOfEdgesToFix}\;
    		    \candidateEnergy $\leftarrow$ \calculateEnergy{\zCand, \wCand} \;
    		    
    		    \If{\candidateEnergy \textless \currentEnergy \KwOr \shouldAcceptCandidate{}}{
    		        \z $\leftarrow$ \zCand, \w $\leftarrow$ \wCand\;
    		    }
    		   
    		    \T $\leftarrow$ \coolingSchedule{\k}\;
    		  
		}
		\Return no complementary tours found;
  }\;
\Proc{\TestVertexAdjacency{\x, \y, \combinedGraph, \startTemp, \numOfEdgesToFix}}{
	\z, \w  $\leftarrow$ \SimulatedAnnealing{\x, \y, \combinedGraph, \startTemp, \numOfEdgesToFix}
	
	\eIf{\z and \w are not empty}{
		\Return vertices $x^v$ and $y^v$ are not adjacent\;
	}{
		\Return vertices $x^v$ and $y^v$  are probably adjacent\;
	}
}
\end{algorithm}

\section{Generation of the initial state} 

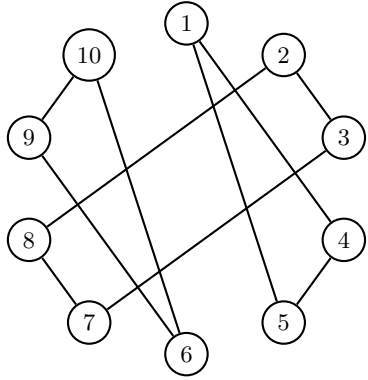
\begin{figure}[h]
	\centering
	\begin{tikzpicture}[scale=0.6]
	\begin{scope}[every node/.style={circle,thick,draw}]
	\graph [clockwise] {
		subgraph I_n [n=10,name=A, radius=2.2cm]; 
		
		A 1 --[thick] A 4;
		A 4 --[thick] A 5; 
		A 5 --[thick] A 1;
		
		A 2 --[thick] A 3;
		A 3 --[thick] A 7;
		A 7 --[thick] A 8;
		A 8 --[thick] A 2;
		
		A 9 --[thick] A 10;
		A 10 --[thick] A 6; 
		A 6 --[thick] A 9;
		
	};
	\end{scope}
	\end{tikzpicture}
	\caption{A vertex-disjoint cycle cover}
	\label{image3}
\end{figure}

To generate the initial system state and neighbor candidate states, we construct a vertex-disjoint cycle cover of the multigraph $x \cup y$ (Fig.~\ref{image3}). 

If $x$ and $y$ are undirected Hamiltonian cycles, then all vertices in the multigraph $x \cup y$ have degrees equal to $4$. Let $z$ be a vertex-disjoint cycle cover of $x \cup y$, then all the remaining edges form a graph $w = (x \cup y) \backslash z$ with all vertex degrees being equal to $2$. Thus, $w$ is also a vertex-disjoint cycle cover of $x \cup y$. 

If $x$ and $y$ are directed Hamiltonian tours, then all vertices in the multigraph $x \cup y$ have both indegrees and outdegrees equal to $2$. Let $z$ be a vertex-disjoint cycle cover of $x \cup y$, then in the digraph $w = (x \cup y) \backslash z$ all vertices have both indegrees and outdegrees equal to $1$. Thus, $w$ is also a vertex-disjoint cycle cover of $x \cup y$. 

Finding a vertex-disjoint cycle cover of both directed and undirected graph can be performed in polynomial time by a reduction to perfect matching \cite{Tutte}.
Let us recall that a perfect matching is a set of pairwise nonadjacent edges which matches all vertices of the graph. 
The procedures for directed and undirected graphs are somewhat different. We consider them separately. 

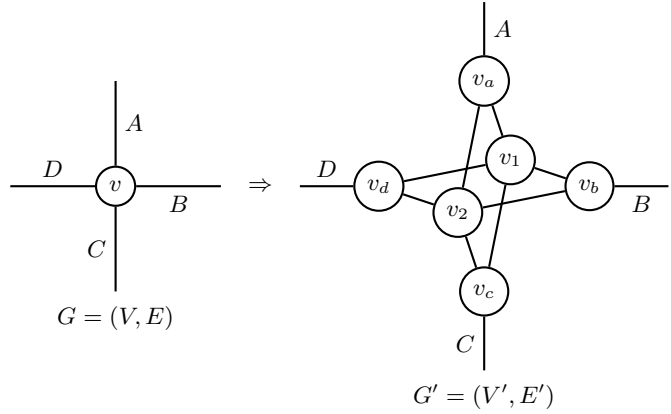
\begin{figure}[t]
	\centering
	\begin{tikzpicture}[scale=0.7]
	\begin{scope}[every node/.style={circle,thick,draw}]
	\node (A) at (0,0) {$v$};
	\end{scope}
	\draw [thick] (A) -- node [right] {$A$} (0,2);
	\draw [thick] (A) -- node [below] {$B$} (2,0);
	\draw [thick] (A) -- node [left] {$C$} (0,-2);
	\draw [thick] (A) -- node [above] {$D$} (-2,0);
	\draw (0, -2.5) node{$G=(V,E)$};
	\node at (2.75,0) {$\Rightarrow$};
	\begin{scope}[xshift=7cm]
	\begin{scope}[every node/.style={circle,thick,draw}]
	\node (A) at (0,2) {$v_a$};
	\node (B) at (2,0) {$v_b$};
	\node (C) at (0,-2) {$v_c$};
	\node (D) at (-2,0) {$v_d$};
	\node (E) at (0.5,0.5) {$v_1$};
	\node (F) at (-0.5,-0.5) {$v_2$};
	\end{scope}
	\draw [thick] (E) edge (A);
	\draw [thick] (E) edge (B);
	\draw [thick] (E) edge (C);
	\draw [thick] (E) edge (D);
	\draw [thick] (F) edge (A);
	\draw [thick] (F) edge (B);
	\draw [thick] (F) edge (C);
	\draw [thick] (F) edge (D);
	\draw [thick] (A) -- node [right] {$A$} (0,3.5);
	\draw [thick] (B) -- node [below] {$B$} (3.5,0);
	\draw [thick] (C) -- node [left] {$C$} (0,-3.5);
	\draw [thick] (D) -- node [above] {$D$} (-3.5,0);
	\draw (0, -4) node{$G'=(V',E')$};
	\end{scope}
	\end{tikzpicture}
	\caption{Construction of the graph $G'$ for the symmetric problem}
	\label{image4}
\end{figure}

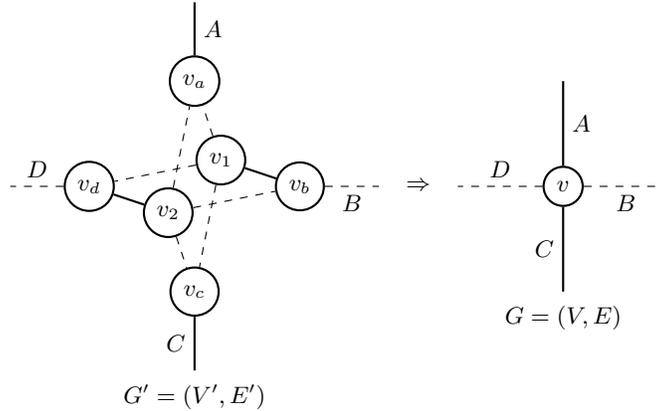
\begin{figure}[t]
	\centering
	\begin{tikzpicture}[scale=0.7]
	\begin{scope}[every node/.style={circle,thick,draw}]
	\node (A) at (0,2) {$v_a$};
	\node (B) at (2,0) {$v_b$};
	\node (C) at (0,-2) {$v_c$};
	\node (D) at (-2,0) {$v_d$};
	\node (E) at (0.5,0.5) {$v_1$};
	\node (F) at (-0.5,-0.5) {$v_2$};
	\end{scope}
	\draw [dashed] (E) edge (A);
	\draw [thick] (E) edge (B);
	\draw [dashed] (E) edge (C);
	\draw [dashed] (E) edge (D);
	\draw [dashed] (F) edge (A);
	\draw [dashed] (F) edge (B);
	\draw [dashed] (F) edge (C);
	\draw [thick] (F) edge (D);
	\draw [thick] (A) -- node [right] {$A$} (0,3.5);
	\draw [dashed] (B) -- node [below] {$B$} (3.5,0);
	\draw [thick] (C) -- node [left] {$C$} (0,-3.5);
	\draw [dashed] (D) -- node [above] {$D$} (-3.5,0);
	\draw (0, -4) node{$G'=(V',E')$};
	
	\node at (4.25,0) {$\Rightarrow$};
	
	\begin{scope}[xshift=7cm]
	\begin{scope}[every node/.style={circle,thick,draw}]
	\node (A) at (0,0) {$v$};
	\end{scope}
	\draw [thick] (A) -- node [right] {$A$} (0,2);
	\draw [dashed] (A) -- node [below] {$B$} (2,0);
	\draw [thick] (A) -- node [left] {$C$} (0,-2);
	\draw [dashed] (A) -- node [above] {$D$} (-2,0);
	\draw (0, -2.5) node{$G=(V,E)$};

	\end{scope}
	\end{tikzpicture}
	\caption{A perfect matching in $G'$ and a vertex-disjoint cycle cover in $G$}
	\label{image5}
\end{figure}

Let $x$ and $y$ be undirected Hamiltonian cycles. 
\begin{enumerate}
    \item [Step 1.] From the multigraph $x \cup y = G = (V, E)$, we construct a new graph $G' = (V', E')$. With each vertex $v \in V$ we associate a gadget $G_v$ that is a complete bipartite subgraph $K_{4,2}$ (note that the degree of $v$ equals 4) as it is shown in Fig.~\ref{image4}:
	\begin{itemize}
	    \item there are $4$ vertices in the outer part ($v_a$, $v_b$, $v_c$ and $v_d$) that correspond to $4$ edges incident to $v$ in $G$ (edges $A$, $B$, $C$, $D$); these vertices are connected with other gadgets; 
		\item there are $2$ vertices in the inner part ($v_1$ and $v_2$) that are connected only with the vertices of the outer part. 
	\end{itemize}
	\item [Step 2.] A perfect matching in $G'$ corresponds to a vertex-disjoint cycle cover in the original graph $G$.
	Indeed, a perfect matching has to cover both inner vertices $v_1$ and $v_2$. 
	Therefore, it includes exactly one edge of \{$(v_1,v_a)$, $(v_1,v_b)$, $(v_1,v_c)$, $(v_1,v_d)$\} and exactly one edge of $\{(v_2,v_a), (v_2,v_b), (v_2,v_c), (v_2,v_d)\}$.
	Both of these edges cover exactly two vertices of $\{v_a,v_b,v_c,v_d\}$.
	The other two vertices has to be covered by the edges that correspond to the edges of $G$ (Fig.~\ref{image5}). 
	We include these edges into $z$, then the degree of each vertex $v$ in the graph $z$ equals $2$, and thus, $z$ is a vertex-disjoint cycle cover of the multigraph $x \cup y$.
\end{enumerate}

A perfect matching in a general undirected graph can be found by Edmond's algorithm \cite{Edmonds} in $O(V^2 E)$ time or using Micali-Vazirani matching algorithm \cite{Micali} in $O(\sqrt{V} E)$ time. We have chosen Edmond's algorithm as a more simple one to implement. 
Note that replacing it with a more efficient Micali-Vazirani algorithm does not require changing the rest of the algorithm.

\begin{figure}[t]
	\centering
	\begin{tikzpicture}[scale=0.75]
	\begin{scope}[every node/.style={circle,thick,draw,inner sep=3pt}]
	\graph [clockwise] {
		subgraph I_n [n=6,name=A, radius=2.2cm]; 
		
		A 1 --[->,>=stealth,thick] A 2;
		A 1 --[->,>=stealth,thick] A 5; 
		A 2 --[->,>=stealth,thick] A 5;
		A 2 --[->,>=stealth,thick] A 6;
		A 3 --[->,>=stealth,thick] A 1;
		A 3 --[->,>=stealth,thick] A 2;
		A 4 --[->,>=stealth,thick] A 1;
		A 4 --[->,>=stealth,thick] A 3;
		A 5 --[->,>=stealth,thick] A 4;
		A 5 --[->,>=stealth,thick] A 6;
		A 6 --[->,>=stealth,thick] A 3;
		A 6 --[->,>=stealth,thick] A 4;
	};
	\end{scope}
	\node at (3.5,0) {$\Rightarrow$};
	
	\begin{scope}[xshift=5cm,yshift=-2.2cm]
	\begin{scope}[every node/.style={circle,thick,draw,inner sep=3pt}]
	\node (A1) at (0,5) {1};
	\node (A2) at (0,4) {2};
	\node (A3) at (0,3) {3};
	\node (A4) at (0,2) {4};
	\node (A5) at (0,1) {5};
	\node (A6) at (0,0) {6};
	\node (B1) at (3,5) {1};
	\node (B2) at (3,4) {2};
	\node (B3) at (3,3) {3};
	\node (B4) at (3,2) {4};
	\node (B5) at (3,1) {5};
	\node (B6) at (3,0) {6};
	\end{scope}
	
	\node at (0,-1) {$L$};
	\node at (3,-1) {$R$};
	
	\draw [thick] (A1) edge (B2);
	\draw [thick] (A1) edge (B5);
	\draw [thick] (A2) edge (B5);
	\draw [thick] (A2) edge (B6);
	\draw [thick] (A3) edge (B1);
	\draw [thick] (A3) edge (B2);
	\draw [thick] (A4) edge (B1);
	\draw [thick] (A4) edge (B3);
	\draw [thick] (A5) edge (B4);
	\draw [thick] (A5) edge (B6);
	\draw [thick] (A6) edge (B3);
	\draw [thick] (A6) edge (B4);
	\end{scope}
	\end{tikzpicture}
	\caption{Construction of the bipatite graph $D'$ for the asymmetric problem}
	\label{image6}
\end{figure}
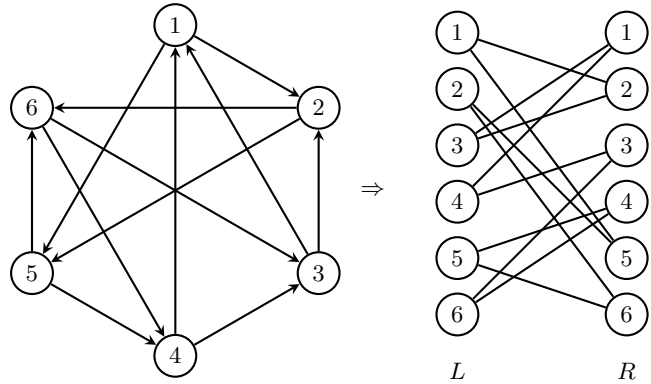

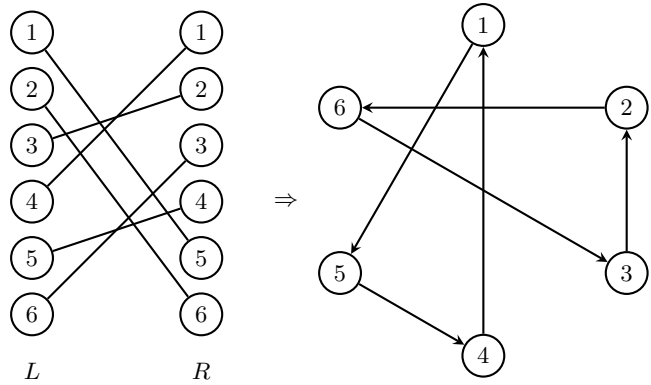
\begin{figure}[t]
	\centering
	\begin{tikzpicture}[scale=0.75]
	\begin{scope}[every node/.style={circle,thick,draw,inner sep=3pt}]
	\node (A1) at (0,5) {1};
	\node (A2) at (0,4) {2};
	\node (A3) at (0,3) {3};
	\node (A4) at (0,2) {4};
	\node (A5) at (0,1) {5};
	\node (A6) at (0,0) {6};
	\node (B1) at (3,5) {1};
	\node (B2) at (3,4) {2};
	\node (B3) at (3,3) {3};
	\node (B4) at (3,2) {4};
	\node (B5) at (3,1) {5};
	\node (B6) at (3,0) {6};
	\end{scope}
	
	\node at (0,-1) {$L$};
	\node at (3,-1) {$R$};
	
	\draw [thick] (A1) edge (B5);
	\draw [thick] (A2) edge (B6);
	\draw [thick] (A3) edge (B2);
	\draw [thick] (A4) edge (B1);
	\draw [thick] (A5) edge (B4);
	\draw [thick] (A6) edge (B3);
	
	\node at (4.5,2) {$\Rightarrow$};
	
	\begin{scope}[xshift=8cm,yshift=2.2cm]
	\begin{scope}[every node/.style={circle,thick,draw,inner sep=3pt}]
	\graph [clockwise] {
		subgraph I_n [n=6,name=A, radius=2.2cm]; 
		
		A 1 --[->,>=stealth,thick] A 5; 
		A 2 --[->,>=stealth,thick] A 6;
		A 3 --[->,>=stealth,thick] A 2;
		A 4 --[->,>=stealth,thick] A 1;
		A 5 --[->,>=stealth,thick] A 4;
		A 6 --[->,>=stealth,thick] A 3;
	};
	\end{scope}
	\end{scope}
	\end{tikzpicture}
	\caption{A perfect matching in $D'$ corresponds to a vertex-disjoint cycle cover of $D$}
	\label{image7}
\end{figure}

Let $x$ and $y$ be directed Hamiltonian tours.
\begin{enumerate}
	\item [Step 1.] From the directed multigraph $x \cup y = D = (V, A)$, we construct a bipartite graph $D' = (L, R, E)$. With each vertex $v \in V$ we associate a pair of vertices $v_L \in L$ and $v_R \in R$, and with each edge $(u,v) \in A$ we associate a new edge $(u_L, v_R)$ in the bipartite graph $D'$ (Fig.~\ref{image6}). 
	\item [Step 2.] A perfect matching in the bipartite graph $D'$ corresponds to a vertex-disjoint directed cycle cover in the original graph $D$.
	Indeed, every vertex of $D$ is a head of exactly one edge and a tail of exactly one edge of a perfect matching in $D'$ (Fig.~\ref{image7}). 
\end{enumerate}

A perfect matching in a bipartite graph can be found by Hopcroft–Karp algorithm \cite{Hopkroft} in $O(\sqrt{V} E)$ time. 

\section{Generation of a neighbor candidate state} 
A process of constructing a neighbor candidate state is shown in Procedure~\ref{transit}. 

\SetAlgorithmName{Algorithm}{Algorithm list}{} 
\SetAlgoCaptionSeparator{.} 
\DontPrintSemicolon
\SetKwProg{Proc}{Procedure}{}{End}
\SetKw{Return}{Return}
\begin{algorithm}[h]
	\label{transit}
	\caption{Constructing a neighbor candidate state}
	\SetKwData{z}{$z$}
	\SetKwData{w}{$w$}
	\SetKwData{numOfEdgesToFix}{$fixEdgesN$}
	\SetKwData{generateNeighbourCandidate}{$GenerateNeighbourCandidate$}
	\SetKwFunction{updateFixedEdges}{$UpdateFixedEdgesQueue$}
	\SetKwFunction{runBlossomAlgorithm}{$RunEdmondsAlgorithm$}
	\SetKwFunction{runHopcroftAlgorithm}{$RunHopcroftKarpAlgorithm$}
	\BlankLine
	
	\Proc{\generateNeighbourCandidate{\z, \w, \numOfEdgesToFix}}{
		\updateFixedEdges{\z, \w, \numOfEdgesToFix}\;

		\eIf{tours \z and \w are directed}{
			\runHopcroftAlgorithm{}\;
		} {
			\runBlossomAlgorithm{}\;
		}  

		\Return \z, \w\;
	}
\end{algorithm}

The algorithm receives as input the current state as the vertex-disjoint cycle covers $z$ and $w$, and the parameter $fixEdgesN$ that set the size of a queue of edges that are fixed in the graph $z$ and the corresponding perfect matching.
When this limit is exceeded, the first edge of the queue is deleted. 

\begin {figure} [t]
\centering
\begin{tikzpicture}[scale=0.9]
\begin{scope}[every node/.style={circle,thick,draw}]
\node (A) at (0,0) {1};
\node (B) at (1,1) {2};
\node (C) at (2.5,1) {3};
\node (D) at (3.5,0) {4};
\node (E) at (2.5,-1) {5};
\node (F) at (1,-1) {6};
\end{scope}
\draw [thick, dashed] (A) edge (B);
\draw [thick] (B) edge (F);
\draw [thick] (F) edge (A);
\draw [thick] (C) edge (D);
\draw [thick, dashed] (D) edge (E);
\draw [thick] (E) edge (C);
\draw (1.75, -1.7) node{$z_{1}$};
\end{tikzpicture}
\hspace*{6mm}
\begin{tikzpicture}[scale=0.9]
\begin{scope}[every node/.style={circle,thick,draw}]
\node (A) at (0,0) {1};
\node (B) at (1,1) {2};
\node (C) at (2.5,1) {3};
\node (D) at (3.5,0) {4};
\node (E) at (2.5,-1) {5};
\node (F) at (1,-1) {6};
\end{scope}
\draw [thick] (A) edge (B);
\draw [thick, dashdotted] (B) edge (C);
\draw [thick] (C) edge (A);
\draw [thick] (D) edge (E);
\draw [thick] (E) edge (F);
\draw [thick] (F) edge (D);
\draw (1.75, -1.7) node{$w_{1}$};
\end{tikzpicture}
\\
\begin{tikzpicture}[scale=0.9]
\begin{scope}[every node/.style={circle,thick,draw}]
\node (A) at (0,0) {1};
\node (B) at (1,1) {2};
\node (C) at (2.5,1) {3};
\node (D) at (3.5,0) {4};
\node (E) at (2.5,-1) {5};
\node (F) at (1,-1) {6};
\end{scope}
\draw [thick, dashed] (A) edge (B);
\draw [thick, dashdotted] (B) edge (C);
\draw [thick] (C) edge (D);
\draw [thick, dashed] (D) edge (E);
\draw [thick] (E) edge (F);
\draw [thick] (F) edge (A);
\draw (1.75, -1.7) node{$z_{2}$};
\end{tikzpicture}
\hspace*{6mm}
\begin{tikzpicture}[scale=0.9]
\begin{scope}[every node/.style={circle,thick,draw}]
\node (A) at (0,0) {1};
\node (B) at (1,1) {2};
\node (C) at (2.5,1) {3};
\node (D) at (3.5,0) {4};
\node (E) at (2.5,-1) {5};
\node (F) at (1,-1) {6};
\end{scope}
\draw [thick] (A) edge (B);
\draw [thick] (B) edge (F);
\draw [thick] (F) edge (D);
\draw [thick] (D) edge (E);
\draw [thick] (E) edge (C);
\draw [thick] (C) edge (A);
\draw (1.75, -1.7) node{$w_{2}$};
\end{tikzpicture}
\caption {Generation of a neighbor candidate state}
\label {image_exchange_rule}
\end{figure}
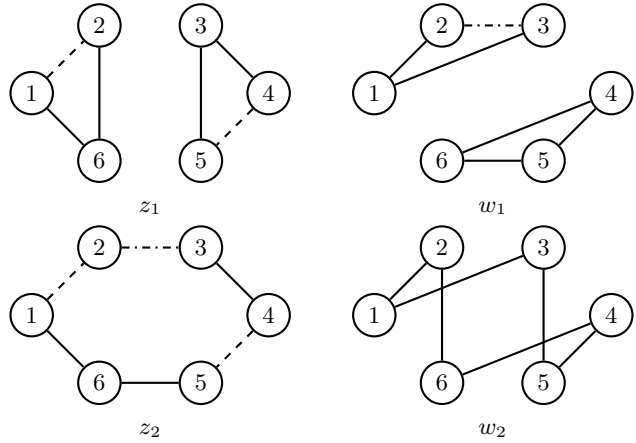

In order to find a neighbor candidate state we chose an edge of $w$ with endpoints in two different connected components of $z$ and add it to the queue of fixed edges (Fig.~\ref{image_exchange_rule}, fixed edges of $z$ are dashed, an edge of $w$ that is added to the queue is dashed dotted).
Such edge always exists due to the connectivity of the multigraph $x \cup y$ . 
If $z$ contains exactly one connected component, then the graphs $z$ and $w$ can be swapped. 
The idea of this procedure is to reduce the number of connected components in $z$ and $w$.
The neighbor candidate state is constructed by the perfect matching algorithms with fixed edges forming the initial matching.

\section{Cooling schedule} 

If a neighbor candidate state has two complementary Hamiltonian tours different from $x$ and $y$ (or any two complementary Hamiltonian tours, if the input is a $2/4$-regular graph), then the algorithm successfully stops and returns a solution.

Otherwise, the energy function is calculated for a neighbor candidate state.
We have chosen the following function
(Procedure~\ref{calcEnergy}): (number of connected components in $z$) + (number of connected components in $w$).

\SetAlgorithmName{Algorithm}{Algorithm list}{} 
\SetAlgoCaptionSeparator{.} 
\DontPrintSemicolon
\SetKwProg{Proc}{Procedure}{}{End}
\SetKw{Return}{Return}
\begin{algorithm}[h]
	\label{calcEnergy}
	\caption{Energy function}
	\SetKwData{z}{$z$}
	\SetKwData{w}{$w$}
	\SetKwData{stateCandidate}{$stateCandidate$}
	\SetKwData{calculateEnergy}{$CalculateEnergy$}
	\SetKwFunction{countComponents}{$CountComponents$}
	
	\BlankLine
	\Proc{\calculateEnergy{\z, \w}}{
		\Return \countComponents{\z} + \countComponents{\w}\;
	}
	\BlankLine
\end{algorithm}

At each step of the algorithm the candidate states correspond to two vertex-disjoint cycle covers $z$ and $w$. Therefore, if the total number of connected components in $z$ and $w$ is equal to $2$, then $z$ and $w$ are Hamiltonian tours. 

If the energy function has decreased compared to the previous state, then we accept a transition to the neighbor candidate state.  

If the energy function has not decreased, then we make a transition with probability
\[P = e^{-\frac{currE - candE}{T}},\]
where $T$ is the current temperature, $currE$ is the current energy value, and $candE$ is the energy of the considered neighbor candidate state.
Such transition is necessary to avoid the problem of falling into the local minimum. 

The current temperature gradually decreases from $initT$ to $0$, and its function depends on the initial temperature $initT$ and the index of the current iteration $k$ (Procedure \ref{tempdecr}).

\SetAlgorithmName{Algorithm}{Algorithm list}{} 
\SetAlgoCaptionSeparator{.} 
\DontPrintSemicolon
\SetKwProg{Proc}{Procedure}{}{End}
\SetKw{Return}{Return}
\begin{algorithm}[!htbp]
	\label{tempdecr}
	\caption{Cooling schedule}
	\SetKwData{k}{$k$}
	\SetKwData{startTemp}{$initT$}
	\SetKwData{coolingSchedule}{$CoolingSchedule$}
	
	\BlankLine
	\Proc{\coolingSchedule{\k}}{
		\Return  \startTemp / \k\;
	}
	\BlankLine
\end{algorithm}



\section{Experiments} 

The algorithm to test vertex adjacencies in the polytopes $\TSP(n)$ and $\ATSP(n)$ was implemented as a console application with different input parameters. Some of them are described below: 

\begin{description}
    \item [\texttt{-{}-}N]~--- number of vertices in the input graph/tour; 
    \item [\texttt{-{}-}times]~--- number of times to run the algorithm; 
    
    \item [\texttt{-{}-}iterN]~--- number of iterations in the simulated annealing algorithm; 
    
    \item [\texttt{-{}-}stateCandidate=random|match]~--- how to generate a neighbor candidate state: 
    \begin{enumerate}
        \item [\textit{random}]: random exchange of edges between tours 
        \item [\textit{match}]: constructing a vertex-disjoint cycle cover and a perfect matching; 
    \end{enumerate}
    
    \item [\texttt{-{}-}exEdgesN]~--- number of edges to randomly exchange between tours (used only for \texttt{-{}-}stateCandidate=random); 
    
    \item [\texttt{-{}-}ansN]~--- multistart: number of repeatedly runs of the algorithm (used only for \texttt{-{}-}stateCandidate=random); 
    
    \item [\texttt{-{}-}fixEdgesN]~--- the size of a queue of edges that can be fixed in the initial matching (used only for \texttt{-{}-}stateCandidate=match). 
\end{description}

We tested the algorithm on random directed and undirected Hamiltonian tours, and also on directed and undirected pyramidal tours.

A Hamiltonian tour 
\[\tau = (1, i_1, i_2, \ldots, i_r, n, j_1, j_2, \ldots, j_{n-r-2})\] 
is called \textit{pyramidal} if 
\[i_1 < i_2 < \ldots < i_r \text{ and } j_1 > j_2 > \ldots > j_{n-r-2}.\]

We chose pyramidal tours for experiments since for them the vertex adjacencies in the corresponding polytopes can be easily verified.
In particular, the following problem: given two pyramidal tours $x$ and $y$, is it possible to construct two complementary tours $z$ and $w$ from the edges of $x$ and $y$, can be solved in linear time \cite{Bondarenko-Nikolaev-2017,Bondarenko-Nikolaev-2018}.
Thus, we run the algorithm on pyramidal tours, for which it is known that the sufficient condition of Lemma~\ref{lemma_nonadjacency} is satisfied.
This allows us to estimate the error percentage when the algorithm could not find complementary tours that are guaranteed to exist. 

The results of the tests for undirected pyramidal tours are presented in Table~\ref{tabl:table1} (Edmond's algorithm is used). 
The algorithm was run with the number of iterations $iterN = 8000$ and the number of fixed edges $fixEdgesN = \lfloor N/3 \rfloor$. 
In the previous version of the program~\cite{Kozlova} a different method to generate neighbor candidate state was implemented~--- the exchange of random edges between two subgraphs $z$ and $w$. Its results are also shown in the table for comparison. For the exchange of random edges the following input parameters were used: the number of iterations $iterN = 50000$, the number of multistart attempts $ansN = 5$ and the number of edges to exchange $exEdgesN = 3$. 

\begin{table}[t]
	\centering
	\caption{Results for undirected pyramidal tours with number of tests $times = 50$}	\label{tabl:table1}
	\begin{tabular}{|*{9}{c|}}
		\hline
		& \multicolumn{4}{c|}{Exchange of random edges} &  \multicolumn{4}{c|}{Reconstructing with perfect matching} \\ 
		\hline
		& Tours & Tours &  &  & Tours & Tours &  & \\
		& are not & are & Average & Accuracy, & are not & are & Average & Accuracy,\\
		& found, & found, & time, & \%, & found, & found, & time, & \%, \\
		N & $TNF_{avg}$ & $TF_{avg}$ & $T_{avg}$ & $Acc$ & $TNF_{avg}$ & $TF_{avg}$ & $T_{avg}$ & $Acc$ \\ 
		\hline
		8 & $-$ & 9,84 & 9,84 & 100 & $-$ & 5,57 & 5,57 & 100 \\ 
		\hline
		16 & 5066,57 & 382,79 & 851,17 & 90 & $-$ & 22,08 & 22,08 & 100 \\ 
		\hline
		24 & 6549,75 & 1403,32 & 5005,82 & 30 & 9096,65 & 33,34 & 214,61 & 98 \\ 
		\hline
		32 & 7832,24 & 1330,37 & 7312,09 & 8 & $-$ & 224,93 & 224,93 & 100 \\ 
		\hline
		40 & 10035,64 & 2351,53 & 9728,28 & 4 & 18656,05 & 610,26 & 971,18 & 98 \\
		\hline
		48 & 13455,39 & $-$ & 13455,39 & 0 & 27695,76 & 978,18 & 3649,94 & 90 \\
		\hline
		64 & 16243,99 & 108,64 & 15921,28 & 2 & 48377,16 & 988,05 & 12361,43 & 76 \\ 
		\hline
		96 &  \multicolumn{4}{c|}{} & 98409,5 & 12293,68 & 53629,27 & 52 \\ 
		\cline{1-1}\cline{6-9}
		128 &  \multicolumn{4}{c|}{\_}  & 158485,17 & 21982,49 & 120264,4 & 28 \\ 
		\cline{1-1}\cline{6-9}
		192 &  \multicolumn{4}{c|}{}  & 334841,38 & 26165,54 & 297800,28 & 12 \\ 
		\hline
	\end{tabular}
\end{table} 

Note that compared to the exchange of random edges, both the accuracy of the algorithm and the size of solved problems have increased. The accuracy can also be adjusted by increasing the number of iterations or changing the maximum number of fixed edges in the queue. 

The results of the tests for directed pyramidal tours and the Hopcroft-Karp algorithm are presented in Table~\ref{tabl:table2}.
The input parameters are similar to the case with undirected tours. 
Here, for almost all considered sizes of test graphs, the algorithm works with the accuracy of $100\%$, which is a very good result. 

\begin{table}[t]
	\centering
	\caption{\label{tabl:table2}Results for directed pyramidal tours (Hopcroft-Karp algorithm), with the number of tests $times = 50$ and number of fixed edges $fixEdgesN = [N/3]$}
	\begin{tabular}{|*{5}{c|}}
		\hline
		& Tours & Tours & Average & \\
		& are not found, & are found, & time, & Accuracy, \%, \\
		N & $TNF_{avg}$ & $TF_{avg}$ & $T_{avg}$ & $Acc$ \\ 
		\hline
		8 & $-$ & 2,27 & 2,27 & 100 \\ 
		\hline
		16 & $-$ & 5,03 & 5,03 & 100 \\ 
		\hline
		24 & $-$ & 19,75 & 19,75 & 100 \\ 
		\hline
		32 & $-$ & 19,14 & 19,14 & 100 \\ 
		\hline
		40 & $-$ & 40,89 & 40,89 & 100 \\ 
		\hline
		48 & $-$ & 95,38 & 95,38 & 100 \\
		\hline
		64 & $-$ & 689,20 & 689,20 & 100 \\
		\hline
		96 & $-$ & 330,21 & 330,21 & 100 \\
		\hline
		128 & 129532,21 & 4514,36 & 9515,07 & 96 \\
		\hline
		192 & 242480,51 & 15783,70 & 70190,93 & 76 \\
		\hline
	\end{tabular}
\end{table}

Finally, Table~\ref{tabl:table3} shows the test results for random directed and undirected Hamiltonian tours. From the table it can be concluded that for random tours the algorithm works even faster than for pyramidal tours.
Note that for undirected graphs the algorithm finds complementary tours more often than for directed graphs. This is due to the fact that $1$-skeleton of the asymmetric traveling salesperson polytope is generally much more dense than $1$-skeleton of the symmetric polytope. For example, the diameter of $1$-skeleton of $\ATSP(n)$ is $2$ \cite{Padberg}, while the best known upper bound for the diameter of $1$-skeleton of $\TSP(n)$ is $4$ \cite{Rispoli}.
Besides, for undirected cycles, the algorithm was able to find a solution for almost all cases. 
We can conclude that for the symmetric traveling salesperson polytope $\TSP(n)$ two random vertices are not adjacent with a very high probability. 


The largest instance that was solved by the algorithm had random Hamiltonian tours on $4096$ vertices and required $1\,530\,682$ ms.
However, due to the long waiting time for several tests, we limited the presented experiments to tours of size under $200$ vertices.

\begin{table}[!htbp]
\centering
\caption{
\label{tabl:table3}Results for random Hamiltonian tours with the number of tests $times = 50$}
\begin{tabular}{|*{9}{c|}}
\hline
& \multicolumn{4}{c|}{Undirected tours} &  \multicolumn{4}{c|}{Directed tours} \\ 
\hline
  & Tours & Tours &  & Percentage & Tours & Tours &  & Percentage \\
  & are not & are & Average & of found & are not & are & Average & of found \\
  & found, & found, & time, & tours, & found, & found, & time, & tours, \\
N & $TNF_{avg}$ & $TF_{avg}$ & $T_{avg}$ & \% & $TNF_{avg}$ & $TF_{avg}$ & $T_{avg}$ & \% \\ 
\hline
8 & 1604,66 & 15,02 & 491,91 & 30 & 3118,99 & 4,05 & 2682,90 & 14 \\ 
\hline
16 & $-$ & 13,76 & 13,756 & 100 & 6799,17 & 4,71 & 4624,95 & 32 \\ 
\hline
24 & $-$ & 28,52 & 28,52 & 100 & 10594,29 & 9,38 & 8265,61 & 22 \\ 
\hline
32 & $-$ & 46,64 & 46,64 & 100 & 14764,29 & 39,19 & 10052,26 & 32 \\ 
\hline
40 & $-$ & 75,27 & 75,27 & 100 & 19184,31 & 19,33 & 15734,62 & 18 \\ 
\hline
48 & $-$ & 87,86 & 87,86 & 100 & 25214,48 & 142,08 & 19197,11 & 24 \\ 
\hline
64 & $-$ & 235,37 & 235,37 & 100 & 38886,56 & 238,95 & 29611,13 & 24 \\ 
\hline
96 & $-$ & 481,37 & 481,37 & 100 & 74654,25 & 1150,20 & 61423,52 & 18 \\ 
\hline
128 & $-$ & 827,42 & 827,42 & 100 & 121790,32 & 1851,29 & 95403,73 & 22 \\ 
\hline
192 & $-$ & 4064,84 & 4064,84 & 100 & 252321,67 & 8979,84 & 213386,98 & 16 \\ 
\hline
\end{tabular}
\end{table}

\section{Conclusion} 

The construction and study of $1$-skeletons of the polytopes associated with intractable problems is of interest for the development and analysis of combinatorial algorithms.
However, for such problems as the traveling salesperson even determining whether two vertices are adjacent or not is an NP-complete problem.
This paper proposes an original heuristic approach based on simulated annealing to verify vertex adjacencies in $1$-skeleton of the traveling salesperson polytope.
The algorithm has a one-sided error: the answer ``not adjacent'' is always correct, while the answer ``probably adjacent'' leaves the possibility that the vertices actually are not adjacent.
The algorithm showed good practical results during the experiments.

\subsubsection*{Acknowledgments.}
The research is supported by the grant of the President of the Russian Federation MK-2620.2018.1. 

%
%
%

\bibliographystyle{splncs04}






\end{document}